\documentclass[twoside]{amsart}

\usepackage{latexsym}
\usepackage{amsthm,amsfonts,amssymb}
\usepackage[leqno]{amsmath}
\usepackage[all]{xy} \SelectTips{eu}{}
\usepackage[hyperfootnotes=false]{hyperref}

\newcommand{\numberseries}{\mdseries}   

\newlength{\thmtopspace}                
\newlength{\thmbotspace}                
\newlength{\thmheadspace}               
\newlength{\thmindent}                  

\newtheoremstyle{bfupright head,slanted body}
                {\thmtopspace}{\thmbotspace}
                {\slshape}{\thmindent}{\bfseries}{.}{\thmheadspace}
                {{\numberseries \thmnumber{(#2) }}\thmnote{#3}}

\newtheoremstyle{bfupright head,upright body}
                {\thmtopspace}{\thmbotspace}
                {\upshape}{\thmindent}{\bfseries}{.}{\thmheadspace}
                {{\numberseries \thmnumber{(#2) }}\thmnote{#3}}

\newtheoremstyle{bfit head,upright body}
                {\thmbotspace}{\thmbotspace}
                {\upshape}{\thmindent}{\upshape}{.}{\thmheadspace}
                {{\numberseries\thmnumber{(#2) }}
                {\bfseries\itshape\thmnote{\negthickspace#3}}}

\newtheoremstyle{fixed bf head,slanted body}
                {\thmtopspace}{\thmbotspace}{\slshape}
                {\thmindent}{\bfseries}{.}{\thmheadspace}
                {{\numberseries \thmnumber{(#2) }}\thmname{#1}\thmnote{ (#3)}}

\newtheoremstyle{fixed bf head,upright body}
                {\thmtopspace}{\thmbotspace}{\upshape}
                {\thmindent}{\bfseries}{.}{\thmheadspace}
                {{\numberseries \thmnumber{(#2)
                    }}\thmname{#1}\thmnote{ (#3)}}

\newtheoremstyle{independent paragraph}
                {\thmtopspace}{\thmbotspace}
                {\upshape}{\parindent}{\upshape}{}{0pt}
                {\thmnote{#3 }}

\newtheoremstyle{subparagraph}
                {\thmbotspace}{\thmbotspace}
                {\upshape}{\parindent}{\upshape}{}{0pt}
                {\thmnote{#3 }}

\newtheoremstyle{notes}
                {\thmtopspace}{\thmbotspace}
                {\ttfamily}{\thmindent}{\ttfamily\small }{}{0pt}
                {\thmnote{#3 }}

\theoremstyle{bfupright head,slanted body}
\newtheorem{res}{}[section]             \newtheorem*{res*}{}
\newtheorem{subres}{}[res]

\theoremstyle{bfupright head,upright body}
\newtheorem{bfhpg}[res]{}               \newtheorem*{bfhpg*}{}

\theoremstyle{bfit head,upright body}
\newtheorem{bfhspg}[subres]{}           \newtheorem*{bfhspg*}{}

\theoremstyle{fixed bf head,slanted body}
\newtheorem{thm}[res]{Theorem}          \newtheorem*{thm*}{Theorem}
\newtheorem{prp}[res]{Proposition}      \newtheorem*{prp*}{Proposition}
        \newtheorem*{cor*}{Corollary}
\newtheorem{lem}[res]{Lemma}            \newtheorem*{lem*}{Lemma}

\theoremstyle{fixed bf head,upright body}
\newtheorem{dfn}[res]{Definition}       \newtheorem*{dfn*}{Definition}
      \newtheorem*{obs*}{Observation}
\newtheorem{rmk}[res]{Remark}           \newtheorem*{rmk*}{Remark}
\newtheorem{exa}[res]{Example}          \newtheorem*{exa*}{Example}

\theoremstyle{independent paragraph}

\theoremstyle{subparagraph}
\newtheorem{spg}{}

\newlength{\thmlistleft}        
\newlength{\thmlistright}       
\newlength{\thmlistpartopsep}   
\newlength{\thmlisttopsep}      
\newlength{\thmlistparsep}      
\newlength{\thmlistitemsep}     

\newcounter{eqc} 
  {\end{list}}%

\newcounter{prt}
\newenvironment{prt}{\begin{list}{\upshape (\alph{prt})}%
    {\usecounter{prt}%
      \setlength{\leftmargin}{\thmlistleft}%
      \setlength{\labelwidth}{\thmlistleft}%
      \setlength{\rightmargin}{\thmlistright}%
      \setlength{\partopsep}{\thmlistpartopsep}%
      \setlength{\topsep}{\thmlisttopsep}%
      \setlength{\parsep}{\thmlistparsep}%
      \setlength{\itemsep}{\thmlistitemsep}}}%
  {\end{list}}%

\newcommand{\prtlbl}[1]{{\upshape(#1)}}

\newcounter{rqm}
  {\end{list}}%

  {\end{list}}%

  {\end{list}}%

\newlength{\cond}
\newenvironment{prf}[1][Proof]{%
\begin{proof}[\bf #1]
\setcounter{equation}{0}
\renewcommand{\theequation}{\arabic{equation}}}
{\end{proof}}

\makeatletter
\def\@nobreak@#1{\mathchoice%
  {\nobreakdef@\displaystyle\f@size{#1}}%
  {\nobreakdef@\nobreakstyle\tf@size{\firstchoice@false #1}}%
  {\nobreakdef@\nobreakstyle\sf@size{\firstchoice@false #1}}%
  {\nobreakdef@\nobreakstyle\ssf@size{\firstchoice@false #1}}%
  \check@mathfonts}%
\def\nobreakdef@#1#2#3{\hbox{{%
                    \everymath{#1}%
                    \let\f@size#2\selectfont%
                    #3}}}%
\makeatother

\renewcommand{\eqref}[1]{\pgref{eq:#1}}

\newcommand{\pgref}[1]{(\ref{#1})}

\newcommand{\exaref}[2][Example~]{#1\pgref{exa:#2}}
\newcommand{\lemref}[2][Lemma~]{#1\pgref{lem:#2}}

\newcommand{\prpref}[2][Proposition~]{#1\pgref{prp:#2}}

\newcommand{\thmref}[2][Theorem~]{#1\pgref{thm:#2}}
\newcommand{\secref}[2][Section~]{#1\ref{sec:#2}}

\newcommand{\partpgref}[2]{(\ref{#1})\prtlbl{#2}}

\newcommand{\partthmref}[3][Theorem~]{#1\partpgref{thm:#2}{#3}}

\newcommand{\corcite}[2][?]{\cite[cor.~#1]{#2}}
\newcommand{\lemcite}[2][?]{\cite[lem.~#1]{#2}}
\newcommand{\prpcite}[2][?]{\cite[prop.~#1]{#2}}

\newcommand{\thmcite}[2][?]{\cite[thm.~#1]{#2}}

\newcommand{\that}[1][\KR]{$#1$-distinguishable}

\newcommand{\Thr}[2]{#2_{\ge #1}}
\newcommand{\Thl}[2]{#2_{\le #1}}

\newcommand{\fb}{\mathfrak{b}}

\newcommand{\Rx}[1][R]{#1/(\x)}
\newcommand{\Ra}[1][R]{#1/(\ba)}
\newcommand{\tp}[3][R]{\nobreak{#2\mspace{-1.5mu}\otimes_{\mspace{-.5mu}#1}\mspace{-1.5mu}#3}}
\newcommand{\stp}[3][R]{\nobreak{#2\mspace{-.5mu}\otimes_{\mspace{-1mu}#1}\mspace{-.25mu}#3}}
\newcommand{\K}{K^{\Rhat}}
\newcommand{\KR}{K^{R}}
\newcommand{\KRh}{K^{R^h}}
\newcommand{\KS}{K^{S}}
\newcommand{\tpK}[2][R]{\tp[#1]{\K\mspace{-2mu}}{#2}}
\newcommand{\stpK}[2][R]{\stp[#1]{\K}{#2}}
\newcommand{\tpKa}[2][R]{\tp[#1]{\K(\ba)\mspace{-2mu}}{#2}}
\newcommand{\tpKx}[2][R]{\tp[#1]{\K(\x)\mspace{-2mu}}{#2}}
\newcommand{\tpKP}[2][R]{\tpP[#1]{\K\mspace{-2mu}}{#2}}
\newcommand{\tpKR}[2][R]{\tp[#1]{\KR\mspace{-2mu}}{#2}}

\newcommand{\tpKS}[2][R]{\tp[#1]{\KS\mspace{-2mu}}{#2}}
\newcommand{\tpKRP}[2][R]{\tpP[#1]{\KR\mspace{-2mu}}{#2}}

\newcommand{\Mor}[2]{\operatorname{Mor}_{\D[\K]}(#1,#2)}

\newcommand{\e}{\varepsilon}
\newcommand{\te}{\hat{\varepsilon}}

\newcommand{\Cat}[2]{{\sf{#2}}(#1)}
\newcommand{\dd}[2]{{\partial}_{#1}^{#2}}
\renewcommand{\H}[2][\no]{\operatorname{H}_{#1}(#2)}

\newcommand{\Susp}[2][]{\mathsf{\Sigma}^{#1}{#2}}
\newcommand{\Hom}[3][R]{\operatorname{Hom}_{#1}(#2,#3)}
\newcommand{\Ext}[4][R]{\operatorname{Ext}_{#1}^{#2}(#3,#4)}

\newcommand{\Tor}[4][R]{\operatorname{Tor}^{#1}_{#2}(#3,#4)}
\newcommand{\D}[1][R]{\Cat{#1}{D}}
\newcommand{\DHom}[3][R]{\operatorname{\mathbf{R}Hom}_{#1}(#2,#3)}
\newcommand{\Dtp}[3][R]{#2\otimes_{#1}^{\mathbf{L}}#3}

\newcommand{\f}{\varphi}

\newcommand{\eq}{\simeq}
\newcommand{\is}{\cong}
\renewcommand{\le}{\leqslant}
\renewcommand{\ge}{\geqslant}
\newcommand{\x}{\pmb{x}}
\newcommand{\ba}{\pmb{a}}
\newcommand{\NNz}{\mathbb{N}_0}
\newcommand{\ZZ}{\mathbb{Z}}
\newcommand{\m}{\mathfrak{m}}

\newcommand{\n}{\mathfrak{n}}
\newcommand{\Rm}{(R,\m)}
\newcommand{\Sn}{(S,\n)}
\newcommand{\Rmk}{(R,\m,k)}

\newcommand{\Rhat}{\widehat{R}}
\newcommand{\Shat}{\widehat{S}}

\newcommand{\lora}{\longrightarrow}
\newcommand{\xla}{\xleftarrow}
\newcommand{\xra}{\xrightarrow}

\newcommand{\xre}{\xra{\;\eq\;}}
\newcommand{\mapdef}[4][\rightarrow]{\mbox{\ensuremath{#2\colon #3 #1 #4}}}
\newcommand{\dmapdef}[4][\lora]{#2\colon #3\:#1\:#4}

\newcommand{\poly}[2][k]{#1[#2]}
\newcommand{\pows}[2][k]{#1[\mspace{-2.3mu}[#2]\mspace{-2.3mu}]}

\newcommand{\Ker}[1]{\mbox{\ensuremath{\operatorname{Ker}#1}}}

\newcommand{\Cone}{\operatorname{Cone}}
\newcommand{\hty}[2][R]{\chi^{#1}_{#2}}
\newcommand{\bid}[2]{\delta^{#1}_{#2}}
\newcommand{\dimR}{\operatorname{dim}R}
\renewcommand{\dim}[2][R]{\operatorname{dim}_{#1}#2}

\newcommand{\rk}[2][k]{\operatorname{rank}_{#1}#2}
\newcommand{\lgt}[2][R]{\operatorname{length}_{#1}#2}
\newcommand{\supP}[1]{\sup{(#1)}}
\newcommand{\infP}[1]{\inf{(#1)}}
\newcommand{\tpP}[3][R]{(\tp[#1]{#2}{#3})}

\newcommand{\lc}[3][\m]{\operatorname{H}^{#2}_{#1}(#3)}

\newcommand{\eqX}{\ensuremath{\mathcal{S}_1}}
\newcommand{\eqY}{\ensuremath{\mathcal{S}_2}}
\newcommand{\eqYY}{\ensuremath{\mathcal{S}_3}}
\newcommand{\eqZ}{\ensuremath{\mathcal{S}_4}}
\newcommand{\cls}[1][C]{\mathsf{#1}}
\newcommand{\MCM}{maximal Cohen--Macaulay }

\renewcommand{\theequation}{\arabic{equation}}
\numberwithin{equation}{res}

\begin{document}

\title{Descent via Koszul extensions}

\dedicatory{Dedicated to Paul C.~Roberts on the occasion of his
  sixtieth birthday}

\author{Lars Winther Christensen}

\address{Department of Math.\ and Stat., Texas Tech University,
  Lubbock, TX 79409, U.S.A.}

\email{lars.w.christensen@math.unl.edu}

\urladdr{http://www.math.ttu.edu/{\tiny $\sim$}lchriste}

\thanks{This work was done while L.W.C.\ visited University of
  Nebraska--Lincoln, partly supported by grants from the Danish
  Natural Science Research Council and the Carlsberg Foundation.}

\author{Sean Sather-Wagstaff}

\address{Department of Mathematical Sciences, Kent State University,
  Mathematics and Computer Science Building, Summit Street, Kent OH
  44242, U.S.A.}

\curraddr{Department of Mathematics, 300 Minard Hall, North Dakota
  State University, Fargo, North Dakota 58105-5075, U.S.A.}

\email{Sean.Sather-Wagstaff@ndsu.edu}

\urladdr{http://math.ndsu.nodak.edu/faculty/ssatherw/}

\thanks{S.S.-W.\ was partially supported by NSF grant NSF~0354281.}

\date{29 February 2008}



\begin{abstract}
  Let $R$ be a commutative noetherian local ring with completion~
  $\Rhat$. We apply differential graded (DG) algebra techniques to
  study descent of modules and complexes from $\Rhat$ to $R'$ where
  $R'$ is either the henselization of $R$ or a pointed \'etale
  neighborhood of $R$: We extend a given $\Rhat$-complex to a DG module
  over a Koszul complex; we describe this DG module equationally and
  apply Artin approximation to descend it to $R'$.
  
  This descent result for Koszul extensions has several applications.
  When $R$ is excellent, we use it to descend the dualizing complex
  from $\Rhat$ to a pointed \'{e}tale neighborhood of $R$; this yields
  a new version of P.~Roberts' theorem on uniform annihilation of
  homology modules of perfect complexes.  As another application we
  prove that the Auslander Condition on uniform vanishing of
  cohomology ascends to $\Rhat$ when $R$ is excellent, henselian, and
  Cohen--Macaulay.
\end{abstract}

\maketitle

\thispagestyle{empty}
\enlargethispage*{\baselineskip}
\section*{Introduction}
\label{sec:intro}

Let $(R,\m)$ be a commutative noetherian local ring with $\m$-adic
completion $\Rhat$. We investigate a pervasive question in local
algebra: When is a given $\Rhat$-module $N$ extended from $R$, i.e.,
when is there an $R$-module $M$ such that $N\is \tp{\Rhat}{M}$?

If $R$ is Cohen--Macaulay, a classical approach to this question is a
two-step analysis that treats the ring and the module separately. Let
$\x$ be a maximal $R$-regular sequence and consider the commutative
diagram of local ring homomorphisms
\begin{equation*}
  \tag{\ensuremath{\ast}}
  \begin{gathered}
    \xymatrix{ R \ar[r] \ar[d] & \Rhat \ar[d]
      \\
      R/(\x) \ar[r]^-{\is} & \Rhat/(\x).  }
  \end{gathered}
\end{equation*}
The bottom map is an isomorphism because $\x$ is a system of
parameters. For every finitely generated $\Rhat$-module $N$, the
module $N/\x N$ is finitely generated over $R/(\x)$ and, hence, also
over $R$. The first step is to identify conditions on $R$ guaranteeing
that $N/\x N$ has the form $M/\x M$ for some finitely generated
$R$-module $M$. The next step is to identify conditions on $N$ such
that the isomorphism $N/\x N \is M/\x M$ forces an isomorphism $N\is
\tp{\Rhat}{M}$. Often, this second step hinges on the good homological
properties of the vertical maps.

If $R$ is not Cohen--Macaulay, then this construction is problematic.
A maximal $R$-regular sequence is not a system of parameters, so the
map $\Rx \to \Rx[\Rhat]$ will not be an isomorphism in general.  This
could be remedied by replacing $\x$ with a system of parameters, but
then the good homological properties of the vertical maps would be
lost. To circumvent these problems, we leave the realm of rings.

We replace the rings $\Rx$ and $\Rx[\Rhat]$ in $(\ast)$ with Koszul
complexes $\KR(\ba)$ and $\K(\ba)$ where $\ba$ is a list of elements
in $\m$ such that $\Ra$ is complete. For example, $\ba$ can be a
system of parameters or a generating sequence for $\m$. The natural
morphism $\KR(\ba) \to \K(\ba)$ is a homology isomorphism and induces
an equivalence between the derived categories $\D[\KR(\ba)]$ and
$\D[\K(\ba)]$ of differential graded (DG) modules. Hence the resulting
commutative diagram of differential graded algebra homomorphisms
supports a two-step analysis parallel to the one described above.

In step one, contained in \thmref[]{approx}, we identify conditions on
$R$ under which a DG module over $\K(\ba)$ that is extended from
$\Rhat$ is also extended from $R$.

\begin{res*}[Theorem~A]
  Let $(R,\m)$ be an excellent henselian local ring and $\ba \in \m$ a
  list of elements such that $\Ra$ is complete.  For every
  $\Rhat$-complex $N$ whose homology is finitely generated over
  $\Rhat$, there is a complex $M$ of~finitely generated $R$-modules
  such that $\tpKa{M}$ is isomorphic to $\tpKa[\Rhat]{N}$ in the
  derived category $\D[\K(\ba)]$.%
\end{res*}

In this theorem, if $N$ is a module then, under additional conditions
on $\ba$ or $N$, also $M$ is a module.  As applications we obtain the
next two theorems. The first is contained in \thmref[]{adsy}; it
extends (the commutative case of) lifting results of Auslander, Ding,
and Solberg \cite{ADS-93}; for definitions see~\pgref{lift}.

\begin{res*}[Theorem~B]
  Let $(R,\m)$ be an excellent henselian local ring and $\x\in\m$ an
  $R$-regular sequence such that $S=\Rx$ is complete.  Let $N$ be a
  finitely generated $S$-module.  If $\,\Ext[S]{2}{N}{N}=0$, then $N$
  has a lifting to $R$.  If $\,\Ext[S]{1}{N}{N}=0$, then any two
  liftings of $N$ to $R$ are isomorphic.
\end{res*}

The second application, contained in \thmref[]{ac}, is an ascent
result for Auslander's conditions on vanishing of cohomology for
finitely generated modules; see \pgref{AC}.

\begin{res*}[Theorem~C]
  Let $R$ be an excellent henselian Cohen--Macaulay local ring. If $R$
  satisfies the (Uniform) Auslander Condition, then so does $\Rhat$.
\end{res*}

In step two of the analysis, we give a condition on $N$ sufficient to
ensure that an isomorphism of DG modules $\tpKa{M} \eq
\tpKa[\Rhat]{N}$ forces an isomorphism of complexes $\tp{\Rhat}{M} \is
N$.  Semidualizing $\Rhat$-complexes \pgref{sdc} satisfy this
condition, and we obtain Theorem~D, which is part of \thmref[]{sdc}.
It subsumes Hinich's result \cite{VHn93} that an excellent henselian
ring admits a dualizing complex; see also Rotthaus \cite{CRt96}.

\begin{res*}[Theorem~D]
  Let $R$ be an excellent henselian local ring. There is a bijective
  correspondence, induced by the functor $\tp{\Rhat}{-}$, between the
  sets of (shift-)isomorphism classes of semidualizing complexes in
  the derived categories $\D$ and $\D[\Rhat]$.%
\end{res*}%

Semidualizing complexes also furnish an example of how the conclusion
of Theorem A may fail for rings that are not excellent and henselian;
see \exaref[]{inadequate}.

\begin{spg}
  Much of this work is done in a setting broader than suggested by the
  above results.  In \thmref[]{sdc} we show that, if $R$ is excellent,
  then every semidualizing $\Rhat$-complex descends to the
  henselization $R^h$ and, moreover, that any finite collection of
  semidualizing $\Rhat$-complexes descends to a pointed \'{e}tale
  neighborhood of $R$.  This allows us to prove, in \thmref[]{PRb}, a
  new version of Roberts' theorem~\cite{PRb76} on uniform annihilation
  of homology modules of perfect complexes.  This, in turn, applies to
  recover a recent result of Zhou \cite{CZh07} on uniform annihilation
  of local cohomology modules.
\end{spg}

As to the organization of the paper, background material is collected
in \secref{algebra}, and Theorem~A is proved in
\secref[Sections~]{approx}--\secref[]{descent}.  Applications,
including Theorems~B and~C, are given in \secref[Sections~]{vanishing}
and \secref[]{Paul}. Theorem~D is proved in \secref{sdc}.

\section{Algebra and module structures}
\label{sec:algebra}

In this paper, $\Rmk$ is a commutative noetherian local ring with
$\m$-adic completion $(\Rhat, \widehat{\m}, k)$.  For a list of
elements $\ba=a_1,\ldots,a_e$ in $\m$, we denote the Koszul complex on
$\ba$ by $\KR(\ba)$.  If $\Ra$ is complete, then we call $\ba$ a
\emph{co-complete sequence}.

For the rest of this section, fix a list of elements $\ba \in\m$ and
set $\KR=\KR(\ba)$.

\begin{bfhpg}[Complexes]
  \label{cx}
  We employ homological grading for complexes of $R$-modules
  \begin{equation*}
    M = \cdots \xra{\dd{n+2}{M}} M_{n+1} \xra{\dd{n+1}{M}} M_n
    \xra{\dd{n}{M}} M_{n-1} \xra{\dd{n-1}{M}} \cdots
  \end{equation*}
  and call them \emph{$R$-complexes} for short. Let $M$ be an
  $R$-complex and $m$ an integer. The \emph{$m$-fold shift} of $M$ is
  denoted $\Susp[m]{M}$; it is the complex with $(\Susp[m]{M})_n =
  M_{n-m}$ and $\partial^{\Susp[m]{M}}_n = (-1)^m\partial^M_{n-m}$.
  The \emph{hard right truncation of $M$ at $m$}, denoted
  $\Thr{m}{M}$, is given by
  \begin{equation*}
    (\Thr{m}{M})_n =
    \begin{cases}
      M_n & \text{if $n \ge m$}\\
      0 & \text{if $n < m$}
    \end{cases}
    \qquad \text{and} \qquad
    \dd{n}{\Thr{m}{M}} =
    \begin{cases}
      \dd{n}{M} & \text{if $n > m$}\\
      0 & \text{if $n\le m$.}
    \end{cases}
  \end{equation*}
  The \emph{hard left truncation of $M$ at $m$} is denoted $M_{\le m}$
  and defined similarly.
  
  A complex $M$ is \emph{bounded} if $M_n =0$ when $|n| \gg 0$.  The
  quantities $\sup{M}$ and $\inf{M}$ are the supremum and infimum,
  respectively, of the set $\{n\in\ZZ \mid \H[n]{M}\ne 0\}$. We say
  that $M$ is \emph{homologically bounded} if $\H{M}$ is bounded, and
  $M$ is \emph{homologically degreewise finite} if each module
  $\H[n]{M}$ is finitely generated. A complex is \emph{homologically
    finite} if it is homologically both bounded and degreewise finite.
  
  Isomorphisms in the category of $R$-complexes are identified by the
  symbol $\is$. Isomorphisms in $\D$, the derived category of the
  category of $R$-modules, are identified by the symbol $\eq$. A
  morphism $\alpha$ between $R$-complexes corresponds to an
  isomorphism in $\D$ if and only if the induced morphism $\H{\alpha}$
  in homology is an isomorphism or, equivalently, the mapping cone
  $\Cone{\alpha}$ is exact; when these conditions are satisfied,
  $\alpha$ is called a \emph{quasiisomorphism}.
  
  Every $R$-complex $M$ has a semifree resolution $P \xre M$; see
  \prpcite[6.6]{fht}.  Such resolutions allow definition of derived
  tensor product and Hom functors $\Dtp{-}{-}$ and $\DHom{-}{-}$
  because the functors $\tp{P}{-}$ and $\Hom{P}{-}$ preserve
  quasiisomorphisms of $R$-complexes.
  
  The Koszul complex $\KR$ is a bounded complex of finite rank free
  $R$-modules, in particular, it is semifree and so the functors
  $\tpKR{-}$ and $\Dtp{\KR}{-}$ are naturally isomorphic.  This fact
  will be used without further mention. If $M$ is homologically
  degreewise finite, then \cite[1.3]{HBFSIn03} provides the
  (in)equalities
  \begin{equation}
    \label{eq:infsup}
    \inf{\tpKR{M}} = \inf{M} \quad\text{and}\quad \sup{M} \le
    \sup{\tpKR{M}} \le e + \sup{M}.
  \end{equation}
  Hence, the complexes $M$ and $\tpKR{M}$ are simultaneously
  homologically bounded.
  
  It is straightforward to verify the following special case of
  tensor-evaluation. For $R$-complexes $M$ and $N$ there is an
  isomorphism in $\D$
  \begin{equation}
    \label{eq:tev}
    \tpKR{\DHom{M}{N}} \xre \DHom{M}{\tpKR{N}}.
  \end{equation}
\end{bfhpg}

\begin{bfhpg}[DG modules over Koszul complexes]
  \label{DGK}
  The Koszul complex $\KR$ can be realized as an exterior algebra, and
  the wedge product endows it with a differential graded
  (DG)\footnote{ Once available, \cite{dga} will be an authoritative
    reference for DG algebra.}  algebra structure that is commutative;
  see e.g.~\prpcite[1.6.2]{bruher}. That is, the product is unitary
  and associative, and it satisfies
  \begin{gather*}
    uv = (-1)^{|u||v|}vu \qquad\text{and}\qquad u^2=0 \text{ when
      $|u|$ is odd}\\
    \dd{}{\KR}(uv) = \dd{}{\KR}(u)v + (-1)^{|u|}u\dd{}{\KR}(v)
  \end{gather*}
  for all $u,v$ in $\KR$, where $|u|$ denotes the degree of $u$.
  
  A \emph{DG $\KR$-module} is an $R$-complex $M$ equipped with a
  $\KR$-multiplication: a morphism of $R$-complexes $\tpKR{M} \to M$,
  written $u\otimes m \mapsto um$, that is unitary and associative and
  satisfies the Leibniz rule
  \begin{equation*}
    \dd{}{M}(um) = \dd{}{\KR}(u)m + (-1)^{|u|}u\dd{}{M}(m)
  \end{equation*}
  for all $u\in\KR$ and $m\in M$. A DG $\KR$-module $M$ is
  \emph{homologically finite} if the homology module $\H{M}$ is
  finitely generated over $\H[0]{\KR}\is \Ra$, equivalently if $M$ is
  homologically finite as an $R$-complex.
  
  If $M$ is an $R$-complex, then $\tpKR{M}$ has a DG $\KR$-module
  structure given by $u(v\otimes x) = (uv)\otimes x$. Moreover, if $M$
  is a homologically finite $R$-complex, then $\tpKR{M}$ is a
  homologically finite DG $\KR$-module.
  
  A morphism of DG $\KR$-modules is a morphism of $R$-complexes that
  is also $\KR$-linear. Isomorphisms in the category of DG
  $\KR$-modules are identified by the symbol $\is$. The derived
  category of the category of DG $\KR$-modules is denoted $\D[\KR]$;
  isomorphisms in this category are identified by the symbol $\eq$. A
  morphism of DG $\KR$-modules corresponds to an isomorphism in
  $\D[\KR]$ if and only if it does so in $\D$ and is then called a
  \emph{quasiisomorphism}.
  
  Every DG $\KR$-module $M$ has a semifree resolution $P \xre M$; see
  \prpcite[6.6]{fht}.  Such resolutions allow definition of derived
  tensor product and Hom functors, $\Dtp[\KR]{-}{-}$ and
  $\DHom[\KR]{-}{-}$ because the functors $\tp[\KR]{P}{-}$ and
  $\Hom[\KR]{P}{-}$ preserve quasiisomorphisms of DG $\KR$-modules.
\end{bfhpg}

\begin{bfhpg}[Local homomorphisms and Koszul complexes]
  Let $\mapdef{\vartheta}{\Rm}{\Sn}$ be a local ring homomorphism,
  that is, $\vartheta(\m) \subseteq \n$. Set
  $\KS=\KS(\vartheta(\ba))$; there is then an isomorphism of
  $S$-complexes and of DG $\KR$-modules
  \begin{equation}
    \label{eq:K}
    \tp{S}{\KR} \is \KS.
  \end{equation}
  Assume $\vartheta$ is flat, i.e.\ it gives $S$ the structure of a
  flat $R$-module, and assume $\Rhat \is \Shat$. If $\ba$ is
  co-complete, then $\vartheta$ induces an isomorphism of rings
  $\Ra\is S/(\vartheta(\ba))$ and, further, a quasiisomorphism of
  $R$-complexes
  \begin{equation}
    \label{eq:KK}
    \KR \xre \KS
  \end{equation}
  which also respects the DG algebra structures. In particular, every
  (homologically finite) DG $\KS$-module is a (homologically finite)
  DG $\KR$-module. Moreover, the functor $\tpKS[\KR]{-}$ is an
  equivalence between $\D[\KR]$ and $\D[\KS]$; it conspires with
  \eqref{KK} to yield an isomorphism in $\D[\KS]$
  \begin{equation}
    \label{eq:KKK}
    \KS \xre \tpKS[\KR]{\KS}.
  \end{equation}
  The homology inverse is the multiplication morphism.
\end{bfhpg}

\section{Equational descriptions of Koszul extensions}
\label{sec:approx}

The next lemma is a crucial step towards Theorem~A from the
introduction.

\begin{lem}
  \label{lem:approx}
  Let $(R,\m)$ be a local ring, $m$ a positive integer, and $P$ a
  complex of finite rank free $\Rhat$-modules such that $P_n=0$ when
  $n<0$ or $n>m$. Fix a co-complete sequence $\ba\in\m$ and set
  $\KR=\KR(\ba)$ and $\K=\K(\ba)$.  There exists a finite system
  $\mathcal{S}$ of polynomial equations with coefficients in $R$ such
  that:
  \begin{prt}
  \item The system $\mathcal{S}$ has a solution in $\Rhat$.
  \item A solution to $\mathcal{S}$ in $R$ yields a complex $A$ of
    finite rank free $R$-modules such that $A_n=0$ when $n<0$ or
    $n>m$, and $\tpK{A} \eq \tpK[\Rhat]{P}$ in $\D[\K]$.
  \end{prt}
\end{lem}

The rest of the section is devoted to the proof of this result; the
argument proceeds in ten steps, the first of which sets up notation.

\begin{bfhpg}[Differentials on $\KR$ and $\K$]
  \label{a0}
  \setcounter{equation}{0}%
  Fix a basis $\e_1,\dots,\e_{2^e}$ for $\KR$ over $R$. For each
  $i=1,\dots,2^e$ set $\te_i = 1\otimes \e_i \in \K$, cf.~\eqref{K}.
  The elements $\te_1,\dots,\te_{2^e}$ form a basis for $\K$ over
  $\Rhat$.  The differential $\dd{n}{\KR}$ is given by a matrix of
  size $\binom{e}{n-1} \times \binom{e}{n}$ with entries in $R$
  \begin{equation*}
    \KR_n \xra{[d_{nij}]} \KR_{n-1}.
  \end{equation*}
  Note that $\dd{n}{\KR} = [d_{nij}] = 0$ when $n<1$ or $n>e$. By
  \eqref{K} the matrix $[d_{nij}]$ also describes the $n$th
  differential on $\K$.
  
  Multiplication on the degree $n$ component of $\KR$ by a basis
  vector $\e_h$ is given by a matrix of size $\binom{e}{n + |\e_h|}
  \times \binom{e}{n}$ with entries in $R$
  \begin{equation*}
    \KR_n \xra{[t^h_{nij}]} \KR_{n+|\e_h|}.
  \end{equation*}
  By \eqref{K} the matrices $[t^h_{nij}]$ also describe multiplication
  by $\te_h$ on $\K$.
\end{bfhpg}

\begin{bfhpg}[A resolution of
  $\K\mspace{-2mu}\otimes_{\Rhat}\mspace{-2mu}P$ over $\KR$]
  \label{a1}
  \setcounter{equation}{0}%
  Without loss of generality, we can assume the complex $P$ is
  minimal; that is, $\partial^P(P) \subseteq \m P$. The DG $\K$-module
  $\tpK[\Rhat]{P}$ is a homologically finite DG $\KR$-module through
  \eqref{KK}. By \cite[prop.~2]{DAp99c} there exists a $\KR$-semifree
  resolution $F \xre \tpK[\Rhat]{P}$, such that the differential of
  $\tp[\KR]{k}{F}$ is $0$ and $F^\natural = \coprod_{i\ge 0}
  \Susp[i]{((\KR)^\natural})^{\beta_i}$ with $\beta_i \in \NNz$.
  (Here $F^\natural$ denotes the graded $R$-module underlying the DG
  $\KR$-module $F$.) Applying the functor $\tpK[\KR]{-}$ yields a
  $\K$-semifree resolution
  \begin{equation*}
    \tpK[\KR]{F} \xre \tpK[\KR]{\tpK[\Rhat]{P}},
  \end{equation*}
  and \eqref{KKK} induces a quasiisomorphism of DG $\K$-modules
  \begin{equation*}
    \tpK[\Rhat]{P} \xre \tpK[\KR]{\tpK[\Rhat]{P}}.
  \end{equation*}
  As $\tpK[\KR]{F}$ is $\K$-semifree, \prpcite[6.4]{fht} provides a
  $\K$-morphism
  \begin{equation}
    \label{eq:f}
    \mapdef[\xre]{\f}{\tpK[\KR]{F}}{\tpK[\Rhat]{P}}.
  \end{equation}
  Since $\f$ is a map between semifree DG $\K$-modules, the induced
  map
  \begin{equation*}
    \dmapdef{\tp[\K]{k}{\f}}{\tp[\K]{k}{\tpKP[\KR]{F}}}
    {\tp[\K]{k}{\tpKP[\Rhat]{P}}}
  \end{equation*}
  is a quasi-isomorphism; see \prpcite[6.7]{fht}.  Each of these
  complexes has zero differential ($P$ is minimal), so
  $\tp[\K]{k}{\f}$ is an isomorphism, hence the underlying graded
  $k$-vector spaces have the same rank. The semifreeness of
  $\tpK[\KR]{F}$ and $\tpK[\Rhat]{P}$ over $\K$ implies
  \begin{equation}
    \label{eq:rk}
    \rk[R]{F}_n = \rk[\Rhat]{\tpKP[\KR]{F}}_n = \rk[\Rhat]{\tpKP[\Rhat]{P}}_n
  \end{equation}
  for all $n$. In particular, $F_n=0$ when $n<0$ or $n>m+e$.
\end{bfhpg}

\begin{bfhpg}[DG structures on $F$ and
  $\K\mspace{-2mu}\otimes_{\KR}\mspace{-2mu}F$]
  \label{a2} \setcounter{equation}{0} For each integer $n$, set $r_n =
  \rk[R]{F_n}$ and fix an $R$-basis for $F_n$. The differential
  $\dd{n}{F}$ is given by an $r_{n-1}\times r_n$ matrix
  \begin{equation*}
    F_n \xra{[u_{nij}]} F_{n-1}
  \end{equation*}
  with entries in $R$. Note that $[u_{nij}] = 0$ when $n<1$ or
  $n>m+e$. There is an isomorphism of DG $\K$-modules
  \begin{equation}
    \label{eq:F}
    \tpK[\KR]{F} \is \tp{\Rhat}{F}
  \end{equation}
  cf.~\eqref{K}. In the $\Rhat$-basis induced by the $R$-basis for
  $F$, the $n$th differential on $\tpK[\KR]{F}$ is also given by the
  matrix $[u_{nij}]$.
  
  For each basis vector $\e_h \in \KR$, cf.~\pgref{a0}, multiplication
  by $\e_h$ on $F$ is given by matrices with entries in $R$
  \begin{equation*}
    F_n \xra{[v^h_{nij}]} F_{n+|\e_h|}.
  \end{equation*}
  By \eqref{F} the matrices $[v^h_{nij}]$ also describe multiplication
  by $\te_h$ on $\tpK[\KR]{F}$.
\end{bfhpg}

\begin{bfhpg}[DG structure on
  $\K\mspace{-2mu}\otimes_{\Rhat}\mspace{-2mu}P$]
  \label{a3} \setcounter{equation}{0} For each integer $n$, fix a
  basis for the free $\Rhat$-module $P_n$ and set $s_n =
  \rk[\Rhat]{P_n}$. The $n$th differential of $P$ is then given by an
  $s_{n-1}\times s_n$ matrix with entries in $\Rhat$
  \begin{equation*}
    P_n \xra{[x_{nij}]} P_{n-1}
  \end{equation*}
  which is zero when $n>m$ or $n<1$. In the basis on $\tpK[\Rhat]{P}$,
  coming from the bases chosen for $\K$ and the modules
  $P_0,\dots,P_m$, the differential
  \begin{equation*}
    \dmapdef{\dd{n}{\stpK[\Rhat]{P}}}{\bigoplus_{p=0}^m
      \tp[\Rhat]{\K_{n-p}}{P_p}}{\bigoplus_{q=0}^m
      \tp[\Rhat]{\K_{n-1-q}}{P_q}}
  \end{equation*}
  is given by a block matrix
  \begin{equation}
    \label{eq:b1}
    \dd{n}{\stpK[\Rhat]{P}} =
    \begin{pmatrix}
      [b_{nij}^{00}] & \cdots & [b_{nij}^{0m}]\\
      \vdots & & \vdots\\
      [b_{nij}^{m0}] & \cdots & [b_{nij}^{mm}]
    \end{pmatrix} = [b_{nij}^{qp}]
  \end{equation}
  where
  \begin{equation}
    \label{eq:b2}
    b^{qp}_{nij} =
    \begin{cases}
      d_{(n-p)ij} & \text{if $p=q$} \\
      (-1)^{n-p}x_{pij} & \text{if $p=q+1$} \\
      0 & \text{otherwise.}
    \end{cases}
  \end{equation}
  Note that $\dd{n}{\stpK[\Rhat]{P}}=0$ when $n < 1$ or $n > m+e$.
  
  For each basis vector $\te_h \in \K$, cf.~\pgref{a0}, multiplication
  by $\te_h$ on $\tpK[\Rhat]{P}$ is given by the formula
  $\te_h(f\otimes g) = (\te f) \otimes g$; see \pgref{DGK}.  In the
  chosen basis for $\tpK[\Rhat]{P}$, multiplication by $\te_h$ on the
  degree $n$ component is given by
  \begin{equation*}
    \dmapdef{\bigoplus_{p=0}^m \tp[\Rhat]{\,[t^h_{(n-p)ij}]}{P_p}}
    {\bigoplus_{p=0}^m \tp[\Rhat]{\K_{n-p}}{P_p}}
    {\bigoplus_{p=0}^m \tp[\Rhat]{\K_{n-p+|\te_h|}}{P_p}}.
  \end{equation*}
  Hence, multiplication by $\te_h$ on $\tpK[\Rhat]{P}$ is given by
  matrices with entries in $R$
  \begin{equation*}
    \tpKP[\Rhat]{P}_n \xra{[w^h_{nij}]} \tpKP[\Rhat]{P}_{n+|\e_h|}.
  \end{equation*}
\end{bfhpg}

\begin{bfhpg}[First set of variables]
  \label{a4} \setcounter{equation}{0} We introduce a finite set of
  variables
  \begin{equation*}
    \{ X_{nij} \mid n=1,\dots,m;\ i=1,\dots,s_{n-1};\
    j=1,\dots,s_{n} \}.
  \end{equation*}
  The equality $\dd{n}{P}\dd{n+1}{P}=0$ says that the elements
  $x_{nij}\in\Rhat$ satisfy a system $\eqX$ of quadratic equations in
  the variables $X_{nij}$ with coefficients 1 and 0, namely the system
  coming from the matrix equations
  \begin{equation}
    \label{eq:X}
    [X_{nij}][X_{(n+1)ij}]=[0]\quad \text{for $n=1,\dots,m-1$}.
  \end{equation}
  
  For later reference, define $[B_{nij}^{qp}]$ for $p,q=0,\dots,m$ and
  $n=1,\dots,m$ to be the block matrix described as in \eqref{b1} and
  \eqref{b2} by
  \begin{equation*}
    B_{nij}^{qp} =
    \begin{cases}
      d_{(n-p)ij} & \text{if $p=q$} \\
      (-1)^{n-p}X_{pij} & \text{if $p=q+1$} \\
      0 & \text{otherwise.}
    \end{cases}
  \end{equation*}
\end{bfhpg}

\begin{bfhpg}[The map $\f$]
  \label{a5} \setcounter{equation}{0} Since $\f$ is a morphism, it
  satisfies the equation
  \begin{equation}
    \label{eq:chain f}
    \f_{n-1}\dd{n}{\stpK[\KR]{F}} - \dd{n}{\stpK[\Rhat]{P}}\f_{n}
    =0
  \end{equation}
  for each $n$. The $\K$-linearity of $\f$ means that there are
  equalities
  \begin{equation*}
    \f(\te_h f) = \te_h\f(f)
  \end{equation*}
  for all $f\in\tpK[\KR]{F}$ and $h = 1,\dots, 2^e$. Each map $\f_n$
  is between free $\Rhat$-modules of rank $r_n$, cf.~\eqref{rk}, so it
  is given by an $r_n\times r_n$ matrix with entries in $\Rhat$
  \begin{equation*}
    \tpKP[\KR]{F}_{n} \xra{[y_{nij}]} \tpKP[\Rhat]{P}_{n}
  \end{equation*}
  which is $0$ when $n<0$ or $n>m+e$.  The $\K$-linearity of $\f$ can,
  therefore, be expressed by commutativity of diagrams
  \begin{equation}
    \label{eq:K f}
    \begin{split}
      \xymatrix@C=6em{\tpKP[\KR]{F}_{n} \ar[r]^-{[y_{nij}]}
        \ar[d]^-{[v^{h}_{nij}]}& \tpKP[\Rhat]{P}_{n}
        \ar[d]_-{[w^{h}_{nij}]}\\
        \tpKP[\KR]{F}_{n + |\e_h|} \ar[r]^-{[y_{(n + |\e_h|)ij}]} &
        \tpKP[\Rhat]{P}_{n + |\e_h|} }
    \end{split}
  \end{equation}
  for $n=0,\dots,m+e$ and $h = 1,\dots, 2^e$.
\end{bfhpg}

\begin{bfhpg}[Second set of variables]
  \label{a6} \setcounter{equation}{0} We introduce another finite set
  of variables
  \begin{equation*}
    \{ Y_{nij} \mid n=0,\dots,m+e;\ i=1,\dots,r_n;\
    j=1,\dots,r_{n} \}.
  \end{equation*}
  By \eqref{chain f} the elements $x_{nij},y_{nij}\in\Rhat$ satisfy a
  system $\eqY$ of equations in the variables $X_{nij}, Y_{nij}$ with
  coefficients in $R$, namely the system coming from the matrix
  equations
  \begin{equation}
    \label{eq:Y}
    [Y_{(n-1)ij}][u_{nij}] - [B_{nij}^{qp}][Y_{nij}] =[0] \quad
    \text{for $n=1,\dots,m+e$}.
  \end{equation}
  By \eqref{K f} the elements $y_{nij}$ satisfy a second system
  $\eqYY$ of equations in $Y_{nij}$ with coefficients in $R$, namely
  those coming from the matrix equations
  \begin{equation}
    \label{eq:YY}
    [Y_{(n + |\e_h|)ij}][v^{h}_{nij}] - [w^{h}_{nij}][Y_{nij}] =[0]
  \end{equation}
  for $h=1,\dots,2^e$ and $n=0,\dots,m+e-|e_h|$.
\end{bfhpg}

\begin{bfhpg}[The mapping cone of $\f$]
  \label{a7} \setcounter{equation}{0} The complex $\Cone{\f}$ consists
  of finite rank free $\Rhat$-modules. In the basis for $\Cone{\f}$
  coming from the bases chosen for $\tpK[\KR]{F}$ and
  $\tpK[\Rhat]{P}$, the $n$th differential is given by the block
  matrix
  \begin{equation*}
    \dd{n}{\Cone{\f}} =
    \begin{pmatrix}
      [b_{nij}^{qp}] & [y_{(n-1)ij}]\\ [0] & -[u_{(n-1)ij}]
    \end{pmatrix}
  \end{equation*}
  which is zero when $n<1$ or $n > m+e+1$.  Since $\f$ is a
  quasiisomorphism, the mapping cone is an exact complex of free
  $\Rhat$-modules and bounded (below). Hence, there exists a homotopy
  between $0$ and the identity on $\Cone{\f}$, i.e.\ a degree $1$
  homomorphism $\sigma$ on $\Cone{\f}$ such that
  \begin{equation}
    \label{eq:cone}
    \sigma_{n-1}\dd{n}{\Cone{\f}} + \dd{n+1}{\Cone{\f}}\sigma_n =
    1^{\Cone{\f}}_n
  \end{equation}
  for every $n$. Each map $\sigma_n$ is given by a matrix of size
  $(r_{n+1} + r_{n}) \times (r_{n}+r_{n-1})$ with entries in $\Rhat$
  \begin{equation*}
    (\Cone{\f})_n \xra{[z_{nij}]} (\Cone{\f})_{n+1}
  \end{equation*}
  which is $0$ when $n<0$ or $n>m+e$.
\end{bfhpg}

\enlargethispage*{\baselineskip}

\begin{bfhpg}[Third set of variables]
  \label{a8} \setcounter{equation}{0} We introduce a third finite set
  of variables
  \begin{equation*}
    \{ Z_{nij} \mid n=0,\dots,m+e;\ i=1,\dots,r_{n+1}+r_{n};\
    j=1,\dots,r_{n} + r_{n-1} \}.
  \end{equation*}
  The equation \eqref{cone} means that the elements $x_{nij}$,
  $y_{nij}$, and $z_{nij}$ satisfy a system of equations in $X_{nij}$,
  $Y_{nij}$, and $Z_{nij}$ with coefficients in $R$, namely the system
  $\eqZ$ coming from the matrix equations
  \begin{equation}
    \label{eq:Z}
    [Z_{(n-1)ij)}]
    \begin{pmatrix}
      [B_{nij}^{qp}] & [Y_{(n-1)ij}]\\ [0] & -[u_{(n-1)ij}]
    \end{pmatrix}
    + 
    \begin{pmatrix}
      [B_{(n+1)ij}^{qp}] & [Y_{nij}]\\ [0] & -[u_{nij}]
    \end{pmatrix}
    [Z_{nij}] = [\delta_{ij}]
  \end{equation}
  for $n=0,\dots,m+e+1$, where $\delta_{ij}$ is the Kronecker delta.
\end{bfhpg}

\begin{bfhpg}[Solutions to $\mathcal{S}$]
  \label{a9}
  \setcounter{equation}{0}%
  By construction, the system $\mathcal{S} = \sqcup_{i=1}^{i=4}
  \mathcal{S}_i$ has a solution in $\Rhat$, namely $x_{nij},\
  y_{nij},\ z_{nij}$; see \eqref{X}, \eqref{Y}, \eqref{YY}, and
  \eqref{Z}. This proves~(a).
  
  For part (b), assume that $\mathcal{S}$ has a solution
  $\tilde{x}_{nij},\ \tilde{y}_{nij},\ \tilde{z}_{nij}$ in $R$. In
  view of the isomorphism $\f$, see \eqref{f}, it suffices to show
  that this yields a complex $A$ of finite rank free $R$-modules such
  that $\tpKR{A} \eq F$ in $\D[\KR]$ and $A_n=0$ when $n<0$ or
  $n>m$. For each $n$, let $A_n$ be a free $R$-module of rank $s_n =
  \rk[\Rhat]{P_n}$ and fix an $R$-basis for $A_n$; note that $A_n=0$
  when $n<0$ or $n>m$. For each $n$ let
  \begin{alignat*}{2}
    \dmapdef{\dd{n}{A}&}{A_n}{A_{n-1}}&
    \quad &\text{be given by}\quad \dd{n}{A}= [\tilde{x}_{nij}]\\
    \dmapdef{\tilde{\f}_n&}{F_n}{\tpKRP{A}_n}&
    &\text{be given by}\quad \tilde{\f}_n = [\tilde{y}_{nij}]\\
    \dmapdef{\tilde{\sigma}_n&}{(\Cone{\tilde{\f}})_n}%
    {(\Cone{\tilde{\f}})_{n+1}}& &\text{be given by}\quad
    \tilde{\sigma}_n = [\tilde{z}_{nij}].
  \end{alignat*}
  Since the elements $\tilde{x}_{nij}$ satisfy \eqX, one has
  $\dd{n}{A}\dd{n+1}{A}=0$; so $A$ is a complex. The elements
  $\tilde{x}_{nij},\ \tilde{y}_{nij}$ satisfy \eqY\ and \eqYY, so the
  map $\tilde{\f}$ is a $\KR$-linear morphism of $R$-complexes.
  Moreover, the elements $\tilde{x}_{nij},\ \tilde{y}_{nij},\
  \tilde{z}_{nij}$ satisfy \eqZ, so the map $\tilde{\sigma}$ is a
  homotopy between $0$ and the identity on $\Cone{\tilde{\f}}$. In
  particular, the cone is exact and, therefore, $\tilde{\f}$ is the
  desired quasiisomorphism. \qed
\end{bfhpg}

\section{Descent of Koszul extensions }
\label{sec:descent}

In this section we accomplish step one of the analysis described in
the introduction. In particular, Theorem~A is a special case of
\partthmref[]{approx}{a}.

\begin{bfhpg}[The approximation property]
  \label{approx}
  The ring $R$ is said to have the \emph{approximation property} if it
  satisfies the following: Given any finite system $\mathcal{S}$ of
  polynomial equations with coefficients in $R$ and any positive
  integer $t$, if $\mathcal{S}$ has a solution in $\Rhat$ then it also
  has a solution in $R$, and the solutions are congruent modulo
  $\widehat{\m}^t$.
  
  By work of D.~Popescu \thmcite[(1.3)]{DPp86}, Rotthaus
  \thmcite[1]{CRt90}, and Spivakovsky \thmcite[11.3]{MSp99}, a local
  ring has the approximation property if and only if it is excellent
  and henselian. For example, every local analytic algebra over a
  perfect field has the approximation property \cite{GScUSt74}; see also
  \cite[(1.19)]{yos}.
\end{bfhpg}

\begin{bfhpg}[Henselization]
  \label{Hensel}
  A pointed \'{e}tale neighborhood of $R$ is a flat local homomorphism
  $R \to R'= \poly[R]{X}_\n/(f)$, where $f$ is a monic
  polynomial whose derivative is a unit in $\poly[R]{X}_\n$, and $\n$
  is a prime ideal lying over $\m$. The set of pointed \'{e}tale
  neighborhoods of $R$ forms a filtered direct system
  $\{R_\lambda\mid\lambda \in \Lambda\}$, and the henselization of $R$
  is the limit $R^h = \varinjlim_\lambda R_\lambda$. The natural map
  $R \to R^h$ is a flat local ring homomorphism, and there is an
  isomorphism $\smash{\widehat{R^h} \is \Rhat}$. See \cite[\S
  18]{egaIV4} and \cite{MRn70}.
  
  Assume $R$ is excellent; by \corcite[(18.7.6)]{egaIV4} and
  \pgref{approx} the henselization $R^h$ then has the approximation
  property.  Let $\mathcal{S}$ be a finite set of polynomial equations
  with coefficients in $R$. If $\mathcal{S}$ has a solution in
  $\Rhat$, then $\mathcal{S}$ has a solution in $R^h$, and it follows
  that there is a pointed \'{e}tale neighborhood $R \to R'$, such that
  $\mathcal{S}$ has a solution in~$R'$.
\end{bfhpg}

\enlargethispage*{\baselineskip}

\begin{thm}
  \label{thm:approx}
  Let $(R,\m)$ be an excellent local ring. Fix a co-complete sequence
  $\ba\in\m$ and set $\K=\K(\ba)$.
  \begin{prt}
  \item For every homologically finite $\Rhat$-complex $N$, there
    exists a homologically finite $R^h$-complex $M$ such that
    $\tpK[R^h]{M}\eq\tpK[\Rhat]{N}$ in $\D[\K]$.
  \item For every list of homologically finite $\Rhat$-complexes
    $N^{(1)},\dots,N^{(t)}$ there exists a pointed \'{e}tale
    neighborhood $R \to R'$ and homologically finite $R'$-complexes
    $M^{(1)},\dots,M^{(t)}$ such that $\tpK[R']{M^{(i)}} \eq
    \tpK[\Rhat]{N^{(i)}}$ in $\D[\K]$.
  \end{prt}
\end{thm}

\begin{prf}
  Let $N$ be any homologically finite $\Rhat$ complex and identify it
  with its minimal semifree resolution; see \cite[prop.~2]{DAp99c}.
  After a shift, we may assume that $N_n =0$ for $n<0$. Set
  $s=\sup{N}$ and $m = s + 2e +1$.  Consider the complex $P =
  \Thl{m}{N}$ and the system $\mathcal{S}$ of equations, whose
  existence and solvability in $\Rhat$ is given by \lemref{approx}. It
  is sufficient to prove the following:
  \begin{bfhspg}[Claim]
    \label{sub}
    \sl If the system $\mathcal{S}$ has a solution in $R$, then there
    exists a homologically finite $R$-complex $M$ such that $\tpK{M}
    \eq \tpK[\Rhat]{N}$ in $\D[\K]$.
  \end{bfhspg}
  \noindent Indeed, $R^h$ has the approximation property, see
  \pgref{approx}, so part (a) follows by applying \pgref{sub} to $R =
  R^h$. A list $N^{(1)},\dots,N^{(t)}$ of homologically finite
  $\Rhat$-complexes also results in a finite system $\mathcal{S}^{(1)}
  \sqcup \cdots \sqcup \mathcal{S}^{(t)}$ of polynomial equations.  As
  noted in \pgref{Hensel}, there is a pointed \'{e}tale neighborhood
  $R \to R'$, such that the compound system and, in particular, each
  subsystem $\mathcal{S}^{(i)}$ has a solution in $R'$. Part (b) now
  follows by applying \pgref{sub} to $R = R'$.

  \begin{bfhspg*}[Proof of {\rm \bf \pgref{sub}}]
    By \lemref{approx} there exists a complex $A$ of finite rank free
    $R$-modules such that $A_n=0$ when $n<0$ or $n>m$, and
    \begin{equation}
      \label{eq:A}
      \tpK{A} \eq \tpK[\Rhat]{\Thl{m}{N}}
    \end{equation}
    in $\D[\K]$. Augment $A$ by an $R$-free resolution of
    $\Ker{\dd{m}{A}}$; this yields a complex $M$ of finite rank free
    $R$-modules with $\sup{M} < m$ and $\Thl{m}{M} \is A$. In
    particular, the isomorphism \eqref{A} can be rewritten as
    \begin{equation}
      \label{eq:B}
      \tpK{\Thl{m}{M}} \eq \tpK[\Rhat]{\Thl{m}{N}}.
    \end{equation}
    Next we show that $\supP{\tpK{M}} < m$. Apply $\tpK[\Rhat]{-}$ to
    the triangle
    \begin{equation*}
      \Thl{m}{N} \xra{\xi^N_m} N \lora \Thr{m+1}{N} \lora \Susp{\Thl{m}{N}}
    \end{equation*}
    and inspect the long exact homology sequence
    \begin{equation*}
      \dots \lora \H[i+1]{\tpK[\Rhat]{\Thr{m+1}{N}}} \lora
      \H[i]{\tpK[\Rhat]{\Thl{m}{N}}} \lora \H[i]{\tpK[\Rhat]{N}} \lora
      \cdots.
    \end{equation*}
    The module $\H[i+1]{\tpK[\Rhat]{\Thr{m+1}{N}}}$ vanishes for $i<m$
    while $\H[i]{\tpK[\Rhat]{N}}$ vanishes for $i> s+e$ by
    \eqref{infsup}. Hence
    \begin{equation}
      \label{eq:C}
      \H[i]{\tpK[\Rhat]{\Thl{m}{N}}}=0 \quad\text{when $s+e < i<m$.}
    \end{equation}
    The isomorphisms $(\tpK{M})_i \is (\tpK{\Thl{m}{M}})_i$ for $i\le
    m$ yield the first isomorphism in the next chain; the second is by
    \eqref{B}, and the vanishing is by \eqref{C}.
    \begin{equation*}
      \H[i]{\tpK{M}} \is \H[i]{\tpK{\Thl{m}{M}}} \is
      \H[i]{\tpK[\Rhat]{\Thl{m}{N}}} = 0 \quad\text{when $s+e < i<m$}.
    \end{equation*}
    Since the modules $\H[i]{M}$ are finitely generated, Nakayama's
    lemma implies that
    \begin{equation*}
      \H[i]{M} = 0 \quad\text{when $s+e < i < m$.}
    \end{equation*}
    As $\sup{M}<m$ it follows that $\sup{M}\le s+e$. Hence
    $\supP{\tpK{M}} \le s+2e <m$.
  
    Next we construct a commutative diagram in the category of DG
    $\K$-modules
    \begin{equation}
      \label{eq:D}
      \begin{split}
        \xymatrix{ \tpK{\Thl{m}{M}} \ar[r]^-{\alpha}_-{\eq}
          \ar[d]_-{\tp[]{\K}{\,\xi^M_m}} &
          \tpK[\Rhat]{\Thl{m}{N}} \ar[d]^-{\tp[]{\K}{\,\xi^N_m}}\\
          \tpK{M} \ar@{.>}[r]^-{\rho} & \tpK[\Rhat]{N}.}
      \end{split}
    \end{equation}
    The top horizontal map exists by \eqref{B} and \prpcite[6.4]{fht}
    as $\tpK{\Thl{m}{M}}$ is $\K$-semifree.  To find a morphism $\rho$
    making the diagram commute, consider the triangle
    \begin{equation*}
      \tpK{\Thl{m}{M}} \lora \tpK{M} \lora \tpK{\Thr{m+1}{M}} \lora
      \Susp{\tpK{\Thl{m}{M}}}
    \end{equation*}
    and apply the functor $\Mor{-}{\tpK[\Rhat]{N}} \is
    \H[0]{\DHom[\K]{-}{\tpK[\Rhat]{N}}}$ to obtain the exact sequence
    of abelian groups
    \begin{align*}
      \Mor{\tpK{M}}{\tpK[\Rhat]{N}} &
      \xra{\mspace{6.5mu}\zeta\mspace{5mu}}
      \Mor{\tpK{\Thl{m}{M}}}{\tpK[\Rhat]{N}}\\ & \lora
      \Mor{\Susp[-1]{\tpK{\Thr{m+1}{M}}}}{\tpK[\Rhat]{N}}.
    \end{align*}
    Note that $\Mor{\Susp[-1]{\tpK{\Thr{m+1}{M}}}}{\tpK[\Rhat]{N}}=0$
    because
    \begin{equation*}
      \infP{\Susp[-1]{\tpK{\Thr{m+1}{M}}}} \ge m
      \qquad\text{and}\qquad
      \supP{\tpK[\Rhat]{N}} \le s+e < m.
    \end{equation*}
    Now $\rho$ can be chosen as any preimage of
    $\tpKP{\xi^N_m}\circ\alpha$ under $\zeta$.  The map $\rho$ is a
    quasiisomorphism, as the vertical maps in \eqref{D} induce
    isomorphisms on homology in degrees less than $m$ and
    $\H[i]{\tpK{M}} =0= \H[i]{\tpK[\Rhat]{N}}$ for $i\ge m$.\qedhere
  \end{bfhspg*}
\end{prf}

\begin{rmk}
  By the existence of minimal semifree resolutions
  \cite[prop.~2]{DAp99c} and the equality of infima \eqref{infsup},
  the complex $M$ in \pgref{sub} can be chosen as a minimal complex of
  finite rank free $R$-modules with $M_n=0$ for $n<\inf{N}$.
\end{rmk}

Recall that a finitely generated $R$-module $M$ is \emph{maximal
  Cohen--Macaulay} if the depth of $M$ equals the Krull dimension of
$R$.

\begin{prp}
  \label{prp:mcm}
  Let $(R,\m)$ be an excellent local ring. Fix a co-complete sequence
  $\ba\in\m$ and set $\K=\K(\ba)$.
  \begin{prt}
  \item For every \MCM $\Rhat$-module $N$, there exists a \MCM
    $R^h$-module $M$ such that $\tpK[R^h]{M} \eq \tpK[\Rhat]{N}$ in
    $\D[\K]$.
  \item For every list of \MCM $\Rhat$-modules $N^{(1)},\dots,N^{(t)}$
    there exists a pointed \'{e}tale neighborhood $R \to R'$ and \MCM
    $R'$-modules $M^{(1)},\dots,M^{(t)}$ such that $\tpK[R']{M^{(i)}}
    \eq \tpK[\Rhat]{N^{(i)}}$ in $\D[\K]$.
  \end{prt}
\end{prp}

\begin{prf}
  (a) By \partthmref{approx}{a} there is a homologically finite
  $R^h$-complex $M$ such that $\tpK[R^h]{M} \eq \tpK[\Rhat]{N}$ in
  $\D[\K]$. It suffices to show that $M$ is (isomorphic in $\D[R^h]$
  to) a \MCM $R^h$-module.
  
  Assume for the moment that the sequence $\ba\in\m$ is a system of
  parameters for $R$.  Then $\ba$ is an $N$-regular sequence, and so
  the depth-sensitivity of $\K$ implies that $\H[i]{\tpK[\Rhat]{N}}
  =0$ for $i>0$. Combining this with \eqref{infsup} it readily follows
  that $M$ is (isomorphic in $\D[R^h]$ to) a module:
  \begin{align*}
    0 &= \inf{\tpKP[\Rhat]{N}} = \inf{\tpKP[R^h]{M}} = \inf{M}\\
    & \le \sup{M} \le \sup{\tpKP[R^h]{M}} = \sup{\tpKP[\Rhat]{N}} = 0.
  \end{align*}
  By depth-sensitivity of the Koszul complex $\smash{\KRh}$, the
  equality
  \begin{equation*}
    0 = \sup{\tpKP[R^h]{M}} = \sup{\tpP[R^h]{K^{R^h}\mspace{-6mu}}{M}}
  \end{equation*}
  implies that $M$ is a \MCM $R^h$-module.
  
  Now consider the general situation, wherein we assume only that
  $\Ra$ is complete.  Let $\x$ be a system of parameters for $R$.
  Applying $\tpKx[\Rhat]{-}$ to the isomorphism $\tpK[R^h]{M} \eq
  \tpK[\Rhat]{N}$ yields the first isomorphism below, while the second
  one uses associativity and commutativity of tensor products.
  \begin{align*}
    \tp[\Rhat]{\K(\x)}{\tpKP[R^h]{M}}& \eq
    \tp[\Rhat]{\K(\x)}{\tpKP[\Rhat]{N}} \\
    \tpK[\Rhat]{\tpP[R^h]{\K(\x)}{M}} &\eq
    \tpK[\Rhat]{\tpP[\Rhat]{\K(\x)}{N}}
  \end{align*}
  The computations from the previous paragraph show that
  $\tp[\Rhat]{\K(\x)}{N}$ is isomorphic in $\D[\Rhat]$ to the
  finite-length module $N/(\x)N$.  The second isomorphism above and
  \thmcite[2.3]{SIn99} now yield the next sequence of equalities
  \begin{equation*}
    \supP{\tpK[\Rhat]{\tpP[R^h]{\K(\x)}{M}}} =
    \supP{\tpK[\Rhat]{\tpP[\Rhat]{\K(\x)}{N}}} = e.
  \end{equation*} 
  The complex $\tp[R^h]{\K(\x)}{M}$ has total homology of finite
  length, and so another application of \thmcite[2.3]{SIn99} yields
  \begin{equation*}
    \supP{\tpK[\Rhat]{\tpP[R^h]{\K(\x)}{M}}} = e +
    \supP{\tp[R^h]{\K(\x)}{M}}.
  \end{equation*} 
  It follows that $\supP{\tp[R^h]{\K(\x)}{M}}=0$, and as in the
  previous argument we conclude that $M$ is (isomorphic in $\D[R^h]$
  to) a \MCM $R^h$-module.

  Part (b) is proved similarly using \partthmref{approx}{b}.
\end{prf}

\section{Applications I: Vanishing of cohomology}
\label{sec:vanishing}

\begin{bfhpg}[Liftings]
  \label{lift}
  Let $\x \in \m$ be an $R$-regular sequence and set $S=\Rx$.  A
  \emph{lifting} of a homologically finite $S$-complex $N$ to $R$ is a
  homologically finite $R$-complex $M$ such that $N\eq \Dtp{S}{M}$ in
  $\D[S]$.  Note that, if $N$ is a module, then a lifting of $N$ to
  $R$ is (isomorphic in $\D$ to) a module $M$ with $N \is \tp{S}{M}$
  and $\Tor{\geq 1}{S}{M}=0$.
  
  In \cite{ADS-93} Auslander, Ding, and Solberg show that vanishing of
  the cohomology modules $\Ext[S]{i}{N}{N}$ for $i=1,2$ guarantees
  existence and uniqueness of a lifting of $N$ to $R$ when $R$ is
  complete.  Yoshino extended these results to complexes in
  \cite{YYs97}.
\end{bfhpg}

The next result uses \thmref{approx} to relax the conditions on the
ring in \cite{ADS-93,YYs97}; it contains Theorem~B from the
introduction.  Note that the assumption that $\x$ is co-complete
yields isomorphisms $\Rx[\Rhat] \is \Rx[R^h]\is\Rx$.

\begin{thm}
  \label{thm:adsy}
  Let $(R,\m)$ be an excellent local ring.  Fix a co-complete
  $R$-regular sequence $\x\in\m$ and set $S=\Rx$.
  \begin{prt}
  \item Every homologically finite $S$-complex $N$ with
    $\Ext[S]{2}{N}{N}=0$ lifts to $R^h$.  If $\,\Ext[S]{1}{N}{N}=0$,
    then any two liftings of $N$ to $R^h$ are isomorphic in $\D[R^h]$.
  \item Let $N^{(1)},\ldots,N^{(t)}$ be homologically finite
    $S$-complexes.  If $\,\Ext[S]{2}{N^{(i)}}{N^{(i)}}=0$ for each
    $i=1,\ldots,t$, then there is a pointed \'{e}tale neighborhood
    $R\to R'$ such that each $N^{(i)}$ lifts to $R'$.  If
    $\,\Ext[S]{1}{N^{(i)}}{N^{(i)}}=0$, then any two liftings of
    $N^{(i)}$ to $R'$ are isomorphic in $\D[R']$.
  \end{prt}
\end{thm}

\begin{prf}
  (a) First, assume $\Ext[S]{2}{N}{N}=0$.  By \lemcite[(3.2)]{YYs97}
  there is a homologically finite $\Rhat$-complex $L$, such that $N\eq
  \Dtp[\Rhat]{S}{L}$ in $\D[S]$.  As $R$ is excellent,
  \partthmref{approx}{a} provides a homologically finite $R^h$-complex
  $M$ such that $\tpKx[R^h]{M} \eq \tpKx[\Rhat]{L}$ in $\D[\K(\x)]$.
  Now the augmentation morphism $\K(\x) \xre S$ yields the second of
  the following isomorphisms in $\D[S]$
  \begin{equation*}
    N\eq \Dtp[\Rhat]{S}{L}\eq \Dtp[R^h]{S}{M}.
  \end{equation*} 
  By the isomorphism $S \is R^h/(\x)$ this shows that $M$ is a lifting
  of $N$ to $R^h$.
  
  Next, assume $\Ext[S]{1}{N}{N}=0$ and let $M$ and $M'$ be liftings
  of $N$ to $R^h$.  It follows that $\Dtp[R^h]{\Rhat}{M}$ and
  $\Dtp[R^h]{\Rhat}{M'}$ are both liftings of $N$ to $\Rhat$, and so
  \lemcite[(3.3)]{YYs97} implies
  $\Dtp[R^h]{\Rhat}{M}\eq\Dtp[R^h]{\Rhat}{M'}$ in $\D[\Rhat]$.  Now
  \lemcite[1.10]{AFrSSW07a} yields the desired isomorphism $M\eq M'$
  in $\D[R^h]$.

  (b) The proof is similar, using \partthmref{approx}{b}.
\end{prf}

\begin{bfhpg}[Auslander Conditions]
  \label{AC}
  The ring $R$ is said to satisfy the \emph{Auslander Condition} if
  for every finitely generated $R$-module $M$ there exists an integer
  $b_M$ such that $\Ext{\gg 0}{M}{X}=0$ implies $\Ext{> b_M}{M}{X}=0$
  for every finitely generated $R$-module $X$. Moreover, $R$ satisfies
  the \emph{Uniform Auslander Condition} if there is an integer $b \ge
  0$ such that $\Ext{\gg 0}{M}{X}=0$ implies $\Ext{> b}{M}{X}=0$ for
  all finitely generated $R$-modules $M$ and $X$. Examples of rings
  that do not satisfy the Auslander Condition were first given by
  Jorgensen and \c{S}ega \cite{DAJLMS04}. For a list of rings that are
  known to satisfy the Auslander Condition, see
  \cite[app.~A]{LWCHHlc}.
\end{bfhpg}

The next result contains Theorem~C from the introduction.

\begin{thm}
  \label{thm:ac}
  Let $R$ be an excellent henselian Cohen--Macaulay local ring. The
  completion $\Rhat$ satisfies the (Uniform) Auslander Condition if
  and only if $R$ satisfies the (Uniform) Auslander Condition.
\end{thm}

\begin{prf}
  It is straightforward to verify the ``only if'' part,
  cf.~\prpcite[(5.5)]{LWCHHlc}.
  
  For the ``if'' part, let $N$ and $Y$ be finitely generated
  $\Rhat$-modules and assume that $\Ext[\Rhat]{\gg 0}{N}{Y}=0$. By a
  standard argument we can without loss of generality assume that $N$
  and $Y$ are \MCM $\Rhat$-modules; see \cite{LWCHHld}.  This
  reduction involves replacing $N$ by a high syzygy and $Y$ by a \MCM
  $\Rhat$-module that approximates it in the sense of \thmcite[A]{MAsROB89}.
  
  Fix a co-complete sequence $\ba\in\m$ and set $\K=\K(\ba)$.  By
  \prpref{mcm} there exist finitely generated $R$-modules $M$ and $X$
  such that $\tpK{M} \eq \tpK[\Rhat]{N}$ and $\tpK{X} \eq
  \tpK[\Rhat]{Y}$ in $\D[\K]$.  This accounts for the last equality in
  the computation below; the first three follow by \eqref{infsup},
  \eqref{tev}, and adjointness.
  \begin{align*}
    \inf{\DHom[\Rhat]{N}{Y}} 
    &= \inf{\tpKP[\Rhat]{\DHom[\Rhat]{N}{Y}}} \\
    &= \inf{\DHom[\Rhat]{N}{\tpK[\Rhat]{Y}}} \\
    &= \inf{\DHom[\K]{\tpK[\Rhat]{N}}{\tpK[\Rhat]{Y}}}\\
    &= \inf{\DHom[\K]{\tpK{M}}{\tpK{X}}}
  \end{align*}
  Combined with a parallel computation starting from $\DHom{M}{X}$,
  this yields the second equality in the next sequence
  \begin{equation}
    \label{eq:ext}
    \begin{aligned}
      \sup\{m\in\ZZ \mid \Ext[\Rhat]{m}{N}{Y}=0\} 
      & = -\inf{\DHom[\Rhat]{N}{Y}} \\
      & = -\inf{\DHom{M}{X}} \\
      & = \sup\{m\in\ZZ \mid \Ext{m}{M}{X}=0\}.
    \end{aligned}
  \end{equation}
  In particular, we have $\Ext{\gg 0}{M}{X}=0$. If $R$ satisfies the
  Auslander Condition, then there is an integer $b=b_M$ such that
  $\Ext{>b}{M}{X}=0$, and \eqref{ext} shows that
  $\Ext[\Rhat]{>b}{N}{Y}=0$. It follows that $\Rhat$ satisfies the
  Auslander Condition.  Ascent of the Uniform Auslander Condition
  follows from the same argument.
\end{prf}


\section{Descent of semidualizing complexes}
\label{sec:sdc}

We now focus on step two of the analysis described in the
introduction, namely, transfer of information from Koszul extensions
to $\Rhat$-complexes.

\begin{dfn}
  \label{dfn:that}
  Fix a co-complete sequence $\ba \in \m$ and set $\KR=\KR(\ba)$.  A
  class $\cls$ of homologically finite $R$-complexes is \emph{\that}
  if it satisfies the following property: Given homologically finite
  $R$-complexes $C$ and $X$, if $C$ is in $\cls$ and $\tpKR{X} \eq
  \tpKR{C}$ in $\D[\KR]$, then $X \eq C$ in $\D$.
\end{dfn}

\begin{lem}
  \label{lem:descent}
  Let $(R,\m)$ be a local ring.  Fix a co-complete sequence $\ba\in\m$
  and set $\K=\K(\ba)$.  Given a \that[\K] class $\cls$ of
  homologically finite $\Rhat$-complexes, if $C$ is in $\cls$ and
  there is a homologically finite $R$-complex $B$ such that $\tpK{B}
  \eq \tpK[\Rhat]{C}$ in $\D[\K]$, then $\tp{\Rhat}{B} \eq C$ in
  $\D[\Rhat]$.
\end{lem}

\begin{prf}
  By the assumption on $\cls$, the claim is immediate from the
  isomorphisms
  \begin{equation*}
    \tpK[\Rhat]{\tpP{\Rhat}{B}} \eq \tpK{B} \eq \tpK[\Rhat]{C}. \qedhere
  \end{equation*}
\end{prf}

The main result of this section concerns semidualizing complexes; the
definition is recalled below. This notion is wide enough to encompass
dualizing complexes in the sense of Grothendieck \cite{rad} and the
relative dualizing complexes of Avramov and Foxby \cite{LLAHBF97}; yet
it is narrow enough to admit a rich theory \cite{LWC01a}. This notion
has been introduced independently by several authors; for example in
Wakamatsu's work on generalized tilting modules \cite{TWk88}.

\begin{bfhpg}[Semidualizing complexes]
  \label{sdc}
  A homologically finite $R$-complex $C$ is \emph{semi\-dualizing} if
  the homothety morphism $\mapdef{\hty{C}}{R}{\DHom{C}{C}}$ is an
  isomorphism in~$\D$.  Further, $C$ is \emph{dualizing} in the sense
  of \cite[V.\S 2]{rad} if it is semidualizing and isomorphic in $\D$
  to a bounded complex of injective modules.
\end{bfhpg}

Theorem~D from the introduction is a special case of part (a) in the
next result.

\begin{thm}
  \label{thm:sdc}
  Let $R$ be an excellent local ring.
  \begin{prt}
  \item For every semidualizing $\Rhat$-complex $C$ there exists a
    semidualizing $R^h$-complex $B$ such that $C \eq
    \tp[R^h]{\Rhat}{B}$. In particular, $R^h$ has a dualizing complex.
  \item For every list of semidualizing $\Rhat$-complexes
    $C^{(1)},\dots,C^{(t)}$ there is a pointed \'{e}tale neighborhood
    $R \to R'$ and semidualizing $R'$-complexes
    $B^{(1)},\dots,B^{(t)}$ such that $C^{(i)} \eq
    \tp[R']{\Rhat}{B^{(i)}}$ for $i=1,\dots,t$. In particular, there
    exists a pointed \'{e}tale neighborhood $R \to R'$ such that $R'$
    admits a dualizing complex.
  \end{prt}
\end{thm}

\begin{prf}
  (a) Let $C$ be a semidualizing $\Rhat$-complex.  Fix a co-complete
  sequence $\ba\in\m$ and set $\K=\K(\ba)$.  By
  \partthmref{approx}{a} there exists a homologically finite
  $R^h$-complex $B$ such that $\tpK[R^h]{B} \eq \tpK[\Rhat]{C}$ in
  $\D[\K]$. The class $\cls$ of semidualizing $\Rhat$-complexes is
  \that[\K] by \lemref{Kinj}, so \lemref{descent} yields an
  isomorphism $\tp[R^h]{\Rhat}{B} \eq C$ in $\D[\Rhat]$, and it
  follows by \pgref{tp} that the complex $B$ is semidualizing for
  $R^h$.  Because $\Rhat$ admits a dualizing complex, this shows that
  $R^h$ also admits a dualizing complex; see
  \thmcite[(5.1)]{LLAHBF92}.
  
  (b) The proof is similar, using \partthmref{approx}{b}.
\end{prf}

The next example demonstrates how badly the conclusion of \thmref{sdc}
(and hence \thmref{approx}) can fail when $R$ is not excellent.

\begin{exa}
  \label{exa:inadequate}
  Let $S_0$ be a field of characteristic zero. For $n>0$ let
  $[V_{nij}]$ be a $2\times 3$ matrix of indeterminants and consider
  the complete normal Cohen--Macaulay local domain
  \begin{equation*}
    S_n = \pows[S_{n-1}]{V_{nij}}\big/I_2(V_{nij}).
  \end{equation*}
  By \corcite[4.9(c)]{SSW07} there are exactly $2^n$ distinct
  shift-isomorphism classes of semidualizing complexes in $\D[S_n]$;
  by \corcite[(3.7)]{LWC01a} each class contains a module.  For each
  $n>0$ there exists, by \thmcite[8]{RCH93}, a Cohen--Macaulay local
  unique factorization domain $R_n$ such that $S_n \is \widehat{R_n}$.
  By \prpcite[3.4]{SSW07} each ring $R_n$ has only the trivial
  semidualizing complex $R_n$ up to shift-isomorphism. Hence, the only
  semidualizing $\widehat{R_n}$-complex that descends to $R_n$ is the
  trivial one $\widehat{R_n}$.
\end{exa}

Also \thmref{adsy} has an application to semidualizing complexes. The
next result extends part of \prpcite[4.2]{AFrSSW07b}.

\begin{prp}
  \label{prp:Rx}
  Let $(R,\m)$ be an excellent henselian local ring and $\x\in\m$ a
  co-complete $R$-regular sequence. There is a bijective
  correspondence, induced by the functor \mbox{$\Dtp{\Rx}{-}$},
  between the sets of (shift-)isomorphism classes of semidualizing
  complexes in $\D$ and $\D[\Rx]$.
\end{prp}

\begin{prf}
  This follows directly from \pgref{tp} and \partthmref{adsy}{a}.
\end{prf}

\section{Applications II: Annihilation of (co)homology}
\label{sec:Paul}

For every local ring $R$, the completion $\Rhat$ has a dualizing
complex. If $R$ is excellent then, by \thmref{sdc}, a dualizing
complex is available much closer to $R$. This is the key to the next
result; if $R$ itself has a dualizing complex, then the conclusion
holds by a result of Roberts \thmcite[1]{PRb76}.

\begin{thm}
  \label{thm:PRb}
  Let $R$ be an excellent local ring of Krull dimension $d$. There
  exists a chain of ideals $\fb_d \subseteq \dots \subseteq \fb_1
  \subseteq \fb_0$ satisfying the following conditions:
  \begin{prt}
  \item For each $i=0,\dots, d$ there is an inequality
    $\dim[]{R/\fb_i} \le i$.
  \item If $F = 0 \to F_r \to \dots \to F_0 \to 0$ is a complex of
    finite rank free $R$-modules with $\lgt{\H{F}} < \infty$, then
    $\fb_i$ annihilates $\H[j]{F}$ for each $j \ge r-i$.
  \end{prt}
\end{thm}

\begin{prf}
  By \partthmref{sdc}{b} there exists a pointed \'{e}tale neighborhood
  $R\to R'$ such that $R'$ admits a dualizing complex.  Note that
  $\dimR'=d$. By \thmcite[1 and preceding prop.\ and def.]{PRb76}
  there exists a chain of ideals $\fb'_d \subseteq \dots \subseteq
  \fb'_1 \subseteq \fb'_0$ in $R'$ such that: ($\text{a}'$)
  $\dim[]{R'/\fb'_i} \le i$, and ($\text{b}'$) if $F' = 0 \to F'_r \to
  \dots \to F'_0 \to 0$ is a complex of finite rank free $R'$-modules
  with $\lgt[R']{\H{F'}} < \infty$, then $\fb'_i \H[j]{F'} =0$ for $j
  \ge r-i$.
  
  For each $i=0,\dots,d$ set $\fb_i = \fb'_i \cap R$. It is not
  difficult to verify that $\dim[]{R/\fb_i} = \dim[]{R'/\fb'_i} \le
  i$. Let $F$ be a complex satisfying the hypothesis of (b). The
  complex $F' = \tp{R'}{F}$ satisfies the hypothesis of ($\text{b}'$).
  Indeed, the isomorphisms $R'/\m R' \is k$ and $\H[j]{F'} \is
  \tp{\H[j]{F}}{R'}$ guarantee $\lgt[R']{\H{F'}} < \infty$.  This
  provides $R$-isomorphisms $\H[j]{F'} \is \H[j]{F}$.  The ideal
  $\fb'_i$ annihilates $\H[j]{F'}$ for $j \ge r-i$ and contains
  $\fb_i$; hence $\fb_i$ annihilates $\H[j]{F}$ for $j \ge r-i$.
\end{prf}

\begin{rmk}
  \label{rmk:nishimura}
  By unpublished examples of Nishimura, an excellent local ring need
  not possess a dualizing complex \cite[exa.~6.1]{JIN}, and a ring
  with a dualizing complex need not be excellent \cite[exa.~4.2]{JIN}.
  In view of this, the hypothesis in \thmref[]{PRb} is neither
  stronger nor weaker than the hypothesis in \thmcite[1]{PRb76}.
\end{rmk}

A classical application of Roberts' theorem \thmcite[1]{PRb76} is to
find uniform annihilators of local cohomology modules $\lc{j}{R}$.
Hochster and Huneke~\cite{MHcCHn90, MHcCHn92}, for instance, do this
when $R$ is an equidimensional local ring that admits a dualizing
complex or is unmixed and excellent.  \thmref{PRb} allows us to drop
the unmixedness condition, thus recovering a recent result of Zhou
\corcite[3.3(ii)]{CZh07}:

\begin{rmk}
  Let $R$ be an equidimensional excellent local ring of Krull
  dimension $d>1$. For each $j$ the local cohomology module
  $\lc{j}{R}$ is a direct limit of homology modules in degree $d-j$ of
  Koszul complexes on powers of a system of parameters. Thus, for each
  $i=0,\dots,d$ the ideal $\fb_i$ annihilates $\lc{j}{R}$ for $j \le
  i$. In particular, $\fb_{d-1}$ annihilates $\lc{j}{R}$ for $j <d$
  and $\dim[]{R/\fb_{d-1}} \le d-1$. If $R$ is equidimensional, then
  $\fb_{d-1}$ is not contained in any minimal prime of $R$.  In the
  terminology of \cite{CZh07}, this means that $\fb_{d-1}$ contains a
  uniform cohomological annihilator for $R$.
\end{rmk}

\appendix
\section*{Appendix. Semidualizing complexes are \that}
\label{sec:obj}
\stepcounter{section}

We start by recalling a few facts about semidualizing complexes; see
\pgref{sdc}.

\begin{bfhpg}[Ascent and descent]
  \label{tp}
  Let $R \to S$ be a local ring homomorphism such that $S$ has finite
  flat dimension as an $R$-module and let $C, C'$ be degreewise
  homologically finite $R$-complexes.  The complex $\Dtp{S}{C}$ is
  $S$-semidualizing if and only if $C$ is $R$-semidualizing.
  Furthermore, if $C$ and $C'$ are semidualizing $R$-complexes such
  that $\Dtp{S}{C} \eq \Dtp{S}{C'}$ in $\D[S]$, then $C\eq C'$ in
  $\D$; see \thmcite[4.5 and 4.9]{AFrSSW07a}.
\end{bfhpg}

The following definition is introduced as a convenience for
\lemref{Kascent}.

\begin{dfn}
  \label{dfn:DGsdm}
  Fix a list of elements $\ba \in\m$ and set $\KR=\KR(\ba)$.  A
  \emph{semidualizing DG $\KR$-module} is a homologically finite DG
  $\KR$-module $M$ such that the homothety morphism
  $\mapdef{\hty[\KR]{M}}{\KR}{\DHom[\KR]{M}{M}}$ is an isomorphism in
  $\D[\KR]$.
\end{dfn}
  
The next lemma shows that the class of semidualizing $R$-complexes is
\that, as defined in \pgref{dfn:that}.

\begin{lem}
  \label{lem:Kascent}
  \label{lem:Kinj}
  Let $(R,\m)$ be a local ring and let $C$ and $C'$ be degreewise
  homologically finite $R$-complexes.  Fix a list of elements $\ba
  \in\m$ and set $\KR=\KR(\ba)$.
  \begin{prt}
  \item The DG $\KR$-module $\tp{\KR}{C}$ is $\KR$-semidualizing if
    and only if $C$ is $R$-semidualizing.
  \item If $C$ and $C'$ are semidualizing and $\tp{\KR}{C} \eq
    \tp{\KR}{C'}$ in $\D[\KR]$, then $C\eq C'$ in $\D$.
  \end{prt}
\end{lem}

\begin{prf}
  For brevity set $K=\KR$. Recall from \cite{LWC01a} that a
  homologically finite $R$-complex $X$ is \emph{$C$-reflexive} if
  $\DHom{X}{C}$ is homologically bounded and the biduality morphism
  $\mapdef{\bid{C}{X}}{X}{\DHom{\DHom{X}{C}}{C}}$ is an isomorphism in
  $\D$.
  
  (a) Recall from \pgref{DGK} that $\tp{K}{C}$ is homologically finite
  over $K$ if and only if $C$ is homologically finite over $R$. In the
  commutative diagram
  \begin{equation*}
    \xymatrix@!C=3em{\tp{K}{R} \ar[d]_-{\stp[]{K}{\hty{C}}} \ar[r]_-{\is} & K
      \ar[rr]^-{\hty[K]{\stp[]{K}{C}}} &
      & \DHom[K]{\tp{K}{C}}{\tp{K}{C}} \ar[d]_-{\eq}\\
      \tp{K}{\DHom{C}{C}} \ar[rrr]^-{\eq}& & & \DHom[R]{C}{\tp{K}{C}}}
  \end{equation*}
  the right-hand vertical isomorphism is by adjointness, and the lower
  horizontal one is tensor-evaluation \eqref{tev}. The diagram shows
  that $\hty[K]{\stp{K}{C}}$ is an isomorphism in $\D[K]$ if and only
  if $\tp{K}{\hty{C}}$ is so. The latter is tantamount to $\hty{C}$
  being an isomorphism in $\D$; to see this apply \eqref{infsup} to
  $\Cone{\tpP{K}{\hty{C}}} \is \tp{K}{\Cone{\hty{C}}}$.
  
  (b) Assume $\tp{K}{C} \eq \tp{K}{C'}$ in $\D[K]$. The first
  isomorphism in the following diagram is the homothety morphism
  \begin{equation*}
    K \xre \DHom[K]{\tp{K}{C'}}{\tp{K}{C}} \xla[\beta]{\;\eq\;}
    \tp{K}{\DHom{C'}{C}}
  \end{equation*}
  while the second is the composition of adjunction and
  tensor-evaluation \eqref{tev}.  It follows that
  $\tp{K}{\DHom{C'}{C}}$ is homologically bounded, and hence so is
  $\DHom{C'}{C}$ by \eqref{infsup}. In the commutative diagram
  \begin{equation*}
    \xymatrixcolsep{2.7em}
    \xymatrixrowsep{5em}
    \xymatrix@!C0{
      \tp{K}{C'} \ar[rrrrrrr]^-{\bid{\stp[]{K}{C}}{\stp[]{K}{C'}}}_-{\eq}
      \ar[d]_-{\stp[]{K}{\bid{C}{C'}}} &&&&&&&
      \DHom[K]{\DHom[K]{\tp{K}{C'}}{\stp[]{K}{C}}}{\stp[]{K}{C}}
      \ar[d]_-{\eq}^-{\DHom[]{\beta}{\tp{K}{C}}} & \\
      *{\phantom{\displaystyle\prod}}
      &\tp{K}{\DHom{\DHom{C'}{C}}{C}} \ar[rrrrrrr]^-{\eq} &&&&&&
      *{\phantom{\displaystyle\prod}} &
      \DHom[K]{\tp{K}{\DHom{C'}{C}}}{\tp{K}{C}} 
    }
  \end{equation*}
  the lower horizontal arrow is the composition of adjointness and
  tensor-evaluation. The diagram shows that $\tp{K}{\bid{C}{C'}}$ is
  an isomorphism in $\D[K]$. An application of \eqref{infsup} to
  $\Cone{\tpP{K}{\bid{C}{C'}}} \is \tp{K}{\Cone{\bid{C}{C'}}}$ implies
  that $\bid{C}{C'}$ is an isomorphism in $\D$. This proves that $C'$
  is $C$-reflexive.  Symmetrically, $C$ is $C'$-reflexive, and by
  \thmcite[5.3]{ATY-05} it then follows that $C$ and $C'$ are
  isomorphic up to shift in $\D$.  From \eqref{infsup} we have $C \eq
  C'$.
\end{prf}


\section*{Acknowledgments}
\label{sec:ack}

We thank Srikanth Iyengar and Roger Wiegand for enlightening
conversations related to this work. We also thank the anonymous referee
for thoughtful comments.


\bibliographystyle{amsplain}

  \newcommand{\arxiv}[2][AC]{\mbox{\href{http://arxiv.org/abs/#2}{\sf arXiv:#2
  [math.#1]}}}
  \newcommand{\oldarxiv}[2][AC]{\mbox{\href{http://arxiv.org/abs/math/#2}{\sf
  arXiv:math/#2
  [math.#1]}}}\providecommand{\MR}[1]{\mbox{\href{http://www.ams.org/mathscine%
t-getitem?mr=#1}{MR#1}}}
  \renewcommand{\MR}[1]{\mbox{\href{http://www.ams.org/mathscinet-getitem?mr=#%
1}{MR#1}}}
\providecommand{\bysame}{\leavevmode\hbox to3em{\hrulefill}\thinspace}
\providecommand{\MR}{\relax\ifhmode\unskip\space\fi MR }
\providecommand{\MRhref}[2]{%
  \href{http://www.ams.org/mathscinet-getitem?mr=#1}{#2}
}
\providecommand{\href}[2]{#2}

\end{document}